\documentclass[10pt]{article}
\usepackage{amsfonts}
\usepackage{amssymb}
\usepackage{amsmath}
\usepackage{mathrsfs}
\begin{document}

\mbox {}
\begin{center}{\noindent{\large\bf {On a problem of Terence Tao}}}
\footnote{Supported by the National Natural Science Foundation of
China( 11471017) and the Natural Science Foundation of HuaiHai
Institute of Technology(KQ10002).}
\end{center}
\begin{center}{ Xue-Gong Sun  }
 \end{center}
\begin{center}{
School of  Sciences,
 HuaiHai Institute of Technology,\\
Lian Yun Gang 222005, P. R. China

 Email: xgsunlyg@163.com}\end{center}
\vskip 5pt {\bf Abstract.} In this paper, we solve a problem of
Terence Tao. We prove that for any $K\geq 2$ and sufficiently large
$N$, the number of primes $p$ between $N$ and $(1+\frac{1}{K})N$
such that $\mid kp+ja^{i}+l\mid$ is composite for all $1\leq a, |j|,
k\leq K$, $1\leq i \leq K\log N$ and $l$ in any set
$L_{N}\subseteq\{-KN,\cdots, KN\}$ of cardinality $K$ with
$ja^{i}+l\neq0$ is at least $C_{K}\frac{N}{\log N}$, where $C_{K}>0$
depending only on $K$.

 \vskip3mm
\centerline{\bf 1. Introduction}

Let $p$ be a prime and $n$ be a nonnegative integer. In 1934,
Romanoff \cite{Romanoff} proved that the set of positive odd
integers which can be expressed in the form $2^n + p$ has a positive
proportion in the set of all positive odd numbers. In 1950, van der
Corput \cite{Crocker} proved that there are a positive proportion
odd integers not of the form $2^n + p$. In the same year, using
covering congruences, Erd\H{o}s \cite{Erdos} proved that there is an
infinite arithmetic progression of positive odd integers each of
which has no representation of the form $2^n + p$. In 1975, Cohen
and Selfridge \cite{CS} proved that there exist infinitely many odd
numbers which are neither the sum nor the difference of a power of
two and a prime power.

Recently, using Selberg's sieve method, Tao \cite{Tao} proved that
for any $K\geq 2$ and sufficiently large $N$, the number of primes
$p$ between $N$ and $(1+\frac{1}{K})N$ such that $\mid kp\pm
ja^{i}\mid$ is composite for all $1\leq a, j, k\leq K$ and $1\leq i
\leq K\log N$, is at least $C_{K}\frac{N}{\log N}$, where $C_{K}$ is
a constant depending only on $K$.

On the other hand, Tao \cite{Tao} posed the following problem:

For any $K\geq 2$ and sufficiently large $N$, the number of primes
$p$ between $N$ and $(1+\frac{1}{K})N$ such that $\mid kp+
ja^{i}+l\mid$ is composite for all $1\leq a, |j|, k\leq K$, $1\leq i
\leq K\log N$ and $l$ in some set $L =L_{N}\subseteq\{-KN,\cdots,
KN\}$ of cardinality at most $K$ is at least $C_{K}\frac{N}{\log
N}$, where $C_{K}$ is a constant depending only on $K$.

Using Tao's idea, in this paper we shall solve the above Tao's
problem. More precisely, we establish

{\bf Theorem 1.} For any $K\geq 2$ and sufficiently large $N$, the
number of primes $p$ between $N$ and $(1+\frac{1}{K})N$ such that
$\mid kp+ ja^{i}+l\mid$ is composite for all $1\leq a, |j|, k\leq
K$, $1\leq i \leq K\log N$ and $l$ in any set
$L_{N}\subseteq\{-KN,\cdots, KN\}$ of cardinality $K$ with
$ja^{i}+l\neq0$ is at least  $C_{K}\frac{N}{\log N}$, where
$C_{K}>0$ depending only on $K$.

{\bf Remark 1.} Let $p(K)$ be a prime with $p(K)>K$. In Theorem 1,
we can take $L_{N}=\{p(K),\cdots,Kp(K)\}$. Moreover, we can take
 $L_{N}=\{K!+1,\cdots,(2K-1)!+1\}$.

{\bf Remark 2.} From Theorem 1,  we know that
 for any $K\geq 2$ and sufficiently large $N$, the
number of primes $p$ between $N$ and $(1+\frac{1}{K})N$ such that
$\mid kp+ ja^{i}+l\mid$ is composite for all $2\leq k\leq K$, $1\leq
a, |j|\leq K$, $1\leq i \leq K\log N$ and $l$ in any set
$L_{N}\subseteq\{-KN,\cdots, KN\}$ of cardinality $K$ is at least
 $C_{K}\frac{N}{\log N}$, where
$C_{K}>0$ depending only on $K$.

\centerline{\bf 2. Proofs}

 In this paper, $p$, $q$, $p_{i,j}$, $q_{i,j}$ are all
primes, and the implied constants in $\ll$, $\gg$ are all absolute.

{\bf Lemma 1 \cite{RL}.} Let $x\geq 2$. Then
 $$\log \log x<\sum\limits_{p\leq
x}\frac{1}{p}<\log\log x+1.$$

{\bf Lemma 2 \cite{RL}.} Let $x\geq 59$. Then
$$\frac{x}{\log x}(1+\frac{1}{2\log x})<\pi(x)<\frac{x}{\log x}(1+\frac{3}{2\log x}).$$

By the Brun's theorem, we get

 {\bf Lemma 3.} There exists a constant $A$ such that

 $$\sum_{p,\,mp+1,\,\,\mbox{primes}}\frac{1}{p}\leq A$$
 for every $1<m\leq M$, where the constant $A$ depending only on $M$.

 {\bf Lemma 4.} For any $2\leq a\leq K$ and $M>12K^{3}$, there exists a set $P_{a}$ of some primes
$ p_{a,t}$ in $[\exp\exp((2^{a-1}-1)(A+1)M),\,\,
\exp\exp((2^{a}-1)(A+1)M))$
 with the following properties:

$(1).$ For any $p_{a,t}$, there exists a prime $q_{p_{a,t}}$ such
that $a^{p_{a,t}}\equiv 1\hskip 0.1mm\pmod {q_{p_{a,t}}}$ and
$q_{p_{a,t}}\geq Mp_{a,t}$.

$(2).$ For each $a$, we have
$\sum\limits_{t}\frac{1}{p_{a,t}}\in[M-3, M]$.

{\bf Proof.} For $a=2$, let $P^{*}_{2}$ be the set of primes in the
interval $[\exp\exp((A+1)M),\,\, \exp\exp(3(A+1)M))$ satisfying that
$mp+1$ is composite for every $1\leq m\leq M$.

By Lemma 1 and Lemma 3, we have

\begin{eqnarray*}&&\sum\limits_{p\in P^{*}_{2}}\frac{1}{p}\geq\sum\limits_{\exp\exp((A+1)M)\leq
p<
\exp\exp(3(A+1)M)}\frac{1}{p}\\
&&-\sum_{m=1}^{m=M}\sum\limits_{\exp\exp((A+1)M)\leq
p<\exp\exp(3(A+1)M),mp+1,\,\,primes}\frac{1}{p}\\
 &&\geq 3(A+1)M-(A+1)M-2-MA\\
&&\geq 2M-2.
\end{eqnarray*}

So, we can find a set $P_{2}$ in $P^{*}_{2}$ with
$\sum\limits_{t}\frac{1}{p_{2,t}}\in[M-3, M]$. Let $q_{p_{2,t}}$ be
the largest prime factor $2^{p_{2,t}}-1$. We know all $q_{p_{2,t}}$
are distinct. By the Fermat's little theorem, we know that $p_{2,t}$
divides $q_{p_{2,t}}-1$. On the other hand, we know that
$mp_{2,t}+1$ is composite for every $1\leq m\leq M$. Thus, we get
$q_{p_{2,t}}\geq Mp_{2,t}$.

Now, suppose that $a>2$ and we have chosen disjoint finite sets of
primes $P_{2}$, $\cdots$, $P_{a-1}$ with the stated properties.

Let $P^{*}_{a}$ be the set of primes in the interval
$[\exp\exp((2^{a-1}-1)(A+1)M),\,\, \exp\exp((2^{a}-1)(A+1)M))$
satisfying that $mp+1$ is composite for every $1\leq m\leq M$.

 Let
$\omega_{a}=\prod\limits_{2\leq
i<a}\prod\limits_{t}(q_{p_{i,t}}-1).$ Note
$q_{p_{i,t}}|i^{p_{i,t}}-1$, by Lemma 2, we get

\begin{eqnarray*}&&\frac{\log(\omega_{a})}{\log a}\\
&&\leq\sum\limits_{2\leq
i<a}\sum\limits_{t}\frac{\log(i^{p_{i,t}}-1)}{\log a}\\
&&\leq\sum\limits_{2\leq
i<a}\sum\limits_{t}p_{i,t}\frac{\log i}{\log a}\\
&&\leq \sum\limits_{2\leq i< a}\sum\limits_{t} p_{i,t}\\
&&\leq \sum\limits_{p\leq \exp\exp((2^{a-1}-1)(A+1)M)}p\\
&&\leq 2\exp(-(2^{a-1}-1)(A+1)M)\exp(2\exp((2^{a-1}-1)(A+1)M)).
\end{eqnarray*}

Thus, we have $\Omega(\omega_{a})\leq
\exp(2\exp((2^{a-1}-1)(A+1)M)).$

Moreover, by Lemma 2, there exists at least
$\exp(2\exp((2^{a-1}-1)(A+1)M))$ primes in the interval $[1,\,
\exp(4\exp((2^{a-1}-1)(A+1)M))].$

So, we get

$$\sum\limits_{ p| \omega_{a}}\frac{1}{p}
\leq \sum\limits_{ p\leq \exp(4\exp((2^{a-1}-1)(A+1)M))}\frac{1}{p}
\leq (2^{a-1}-1)(A+1)M+\log4 +1.
$$

 By Lemma 1 and Lemma 3, we have

\begin{eqnarray*}&&\sum\limits_{p\in P^{*}_{a}, \, p\nmid \omega_{a}}\frac{1}{p}\\
&&\geq\sum\limits_{\exp\exp((2^{a-1}-1)(A+1)M)\leq p<
\exp\exp((2^{a}-1)(A+1)M)}\frac{1}{p}\\
&&-\sum_{m=1}^{m=M}\sum\limits_{\exp\exp((2^{a-1}-1)(A+1)M)\leq
p<\exp\exp((2^{a}-1)(A+1)M),mp+1,\,\,primes}\frac{1}{p}\\
&&-\sum\limits_{ p| \omega_{a}}\frac{1}{p}\\
 &&\geq (2^{a}-1)(A+1)M-1-(2^{a-1}-1)(A+1)M-1-MA-(2^{a-1}-1)(A+1)M-\log4 -1\\
&&\geq M-3.
\end{eqnarray*}

Since $p_{a,t}|q_{p_{a,t}}-1$ and $p_{a,t}\nmid \omega_{a}$, we know
all these $q_{p_{a,t}}$ are distinct.

Similar to $a=2$, we can choose a set $P_{a}$ in $P^{*}_{a}$ with
the stated properties.

This completes the proof of Lemma 4.

 {\bf Proof of Theorem 1.} Let

 $$R=\{(j, k, l): 1\leq |j|, k\leq K,\, l\in L_{N}\}$$
 and take $p$ and $q_{p}$ as in Lemma
 4.

By $\sum\limits_{p\in P_{a}}\frac{1}{p}\in[M-3, M]$, we may
partition
 $P_{a}=\bigcup_{j,\,k,\,l}P_{a,\,j,\,k,\,l}$ in such a way that

 $$\sum\limits_{p\in P_{a,\,j,\,k,\,l}}\frac{1}{p}\in[\frac{M}{4K^{3}},\,\frac{M}{3K^{3}}].$$

 Let $W$ be  the quantity
 $W=\prod\limits_{p}q_{p}$. For $p \in P_{a,j,k,l}$, let $I(a,j,k,l)$
 be the smallest integer $i\geq0$ such that $ja^{i}+l\not\equiv0\hskip 1mm\pmod {q_{p}}$,
 we know that $I(a,j,k,l)=0, 1$.

By the Chinese remainder theorem, we can take $(b,W)=1$ satisfying

 $$
 kb+ja^{I(a,j,k,l)}+l\equiv0\hskip -3mm\pmod {q_{p}}
$$
for every $p\in P_{a,j,k,l}$, $2\leq a\leq K$, and $(j, k, l)\in R$.

Let
\begin{eqnarray*}&&Q=\#\{N\leq  m\leq (1+K^{-1})N: m\equiv b\hskip 0.1mm\pmod {W},
m \, \mbox{prime}, \,\mbox{but}\, |km+ja^{i}+l|  \, \mbox{composite for all}\, \\
&& 1\leq i<K\log N, 1\leq a\leq K, (j, k, l)\in R\},\\
&&Q_{N}=\#\{N\leq  m\leq (1+K^{-1})N: m\equiv b\hskip 0.1mm\pmod
{W}, m \, \mbox{prime}\},
\end{eqnarray*}

and$$ Q_{N,a, i, j, k, l}=\#\{N\leq  m\leq (1+K^{-1})N: m\equiv
b\hskip 0.1mm\pmod {W}, m, \, |km+ja^{i}+l|,\, \mbox{primes}\}.
$$

Similar to the proof of Theorem 1.2 \cite{Tao}, we get

$$Q\geq Q_{N}-\sum\limits_{a=2}^{K}\sum\limits_{1\leq i<K\log
N}\sum\limits_{(j,k,l)\in R}Q_{N,a,i,j,k,l}-\sum\limits_{(j,k,l)\in
R}Q_{N, 1,1,j,k,l}-O(\log N).$$

From the prime number theory in arithmetic progressions, we
have

$$Q_{N}\geq c_{1}\frac{N}{W\log N}\prod\limits_{q|W}(1-\frac{1}{q}).$$

Let $P^{*}=\{p: K<p<N^{\frac{1}{8}},(p,W)=1\}$.

By the Selberg's sieve method, we have

\begin{eqnarray*}&&Q_{N,a, i, j, k, l}=\#\{N\leq  m\leq (1+K^{-1})N: m\equiv b\hskip
0.1mm\pmod {W}, m, \, |km+ja^{i}+l|,\, \mbox{primes}\}\\
&&\leq  \#\{N\leq  m\leq (1+K^{-1})N: m\equiv b\hskip 0.1mm\pmod
{W}, m>N^{\frac{1}{8}}, \, |km+ja^{i}+l|>N^{\frac{1}{8}},\,
\mbox{primes}\}+2N^{\frac{1}{8}}\\
&&\leq \#\{1\leq  r\leq (1+K^{-1})\frac{N}{W}+1:
Wr+b>N^{\frac{1}{8}}, \, |kWr+kb+ja^{i}+l|>N^{\frac{1}{8}}, \,
\mbox{primes}\}+2N^{\frac{1}{8}}\\
 &&\leq  \#\{1\leq  r\leq (1+K^{-1})\frac{N}{W}+1:
((Wr+b)(kWr+kb+ja^{i}+l),p)=1, \,
p\in P^{*}\}+2N^{\frac{1}{8}}\\
&&\ll\frac{N}{W}\prod\limits_{p\nmid ja^{i}+l,\,p\in
P^{*}}(1-\frac{2}{p})\prod\limits_{p| ja^{i}+l,\,p\in
P^{*}}(1-\frac{1}{p})\\
&&\ll\frac{N}{W}\prod\limits_{p\in
P^{*}}(1-\frac{2}{p})\prod\limits_{p| ja^{i}+l,\,p\in
P^{*}}(1-\frac{1}{p})(1-\frac{2}{p})^{-1}\\
&&\ll\frac{N}{W\log^{2}N}\prod\limits_{3\leq p\leq
K}(1-\frac{2}{p})^{-1}\prod\limits_{K<q|W}(1-\frac{2}{q})^{-1}\prod\limits_{K<p|
ja^{i}+l,\,p\nmid W}(1+\frac{1}{p})\\
&&\ll\frac{N}{W\log^{2}N}\prod\limits_{3\leq p\leq
K}(1-\frac{2}{p})^{-1}\prod\limits_{K<q|W}(1-\frac{2}{q\mathbf{}})^{-1}\prod\limits_{K<p|
ja^{i}+l}(1+\frac{1}{p}).
\end{eqnarray*}

Now suppose that $2\leq a\leq K$, $(j,k,l)\in R$.

 Note that if $i\equiv I(a,j,k,l)\hskip 0.1mm\pmod {p}$ for some $p\in
P_{a,j,k,l}$, then $q_{p}|km+ja^{i}+l$, so $q_{p}=|km+ja^{i}+l|$.

Thus, we have

 $$\sum\limits_{1\leq i\leq K\log N,i\equiv I(a,j,k,l)\pmod {p}\,\,\mbox{for some}\,\,p\in P_{a,j,k,l}} Q_{N,a, i, j, k, l}
\ll\log N$$.

 Let
$e_{a,j,l}(d)$ denote the smallest positive integer $i$ such that
$ja^{i}+l\equiv 0\hskip 0.1mm\pmod {d}$. Since $p|d\Rightarrow K<p$,
we know that $d|ja^{i}+l$ if and only if $e_{a,j,l}(d)|i$.

Let $$E(x)=\sum\limits_{0<k\leq x}\sum\limits_{\mu^{2}(d)=1,
e_{a,j,k,l}(d)=k, p|d\Rightarrow
K<p}\frac{2^{\omega(d)}}{d}$$.

Similar to the proof of Lemma 7.8 \cite{MBN}, we have

$$E(x)\ll \log ^{2}x.$$

By partial summation, we have $$\sum\limits_{0<k\leq
x}\frac{1}{k}\sum\limits_{\mu^{2}(d)=1,e_{a,j,k,l}(d)=k,
p|d\Rightarrow K<p}\frac{2^{\omega(d)}}{d}\ll 1.$$

So, we get

$$\sum\limits_{\mu^{2}(d)=1, p|d\Rightarrow
K<p}\frac{2^{\omega(d)}}{de_{a,j,l,k}(d)}=\sum\limits_{0<k}\frac{1}{k}\sum\limits_{\mu^{2}(d)=1,e_{a,j,k,l}(d)=k,
p|d\Rightarrow K<p}\frac{2^{\omega(d)}}{d}\ll 1.$$

By the Cauchy-Schwarz inequality and the Selberg's sieve method, we
know

\begin{eqnarray*}&&\sum\limits_{1\leq i<K\log N, \,p\in P_{a,j,k,l}\Rightarrow p\nmid (i-I(a,j,k,l))}
\prod\limits_{K<p|
ja^{i}+l}(1+\frac{1}{p})\\
&&\ll(\sum\limits_{1\leq i<K\log N, \,p\in P_{a,j,k,l}\Rightarrow
p\nmid (i-I(a,j,k,l))}1)^{\frac{1}{2}} (\sum\limits_{1\leq i<K\log
N}\prod\limits_{K<p|
ja^{i}+l}(1+\frac{1}{p})^{2})^{\frac{1}{2}}\\
&&\ll(K\log N\prod\limits_{p\in
P_{a,j,k,l}}(1-\frac{1}{p}))^{\frac{1}{2}}(\sum\limits_{1\leq
i<K\log N}\prod\limits_{K<p|
ja^{i}+l}\frac{2}{p})^{\frac{1}{2}}\\
&&\ll(K\log N\prod\limits_{p\in
P_{a,j,k,l}}(1-\frac{1}{p}))^{\frac{1}{2}}(\sum\limits_{1\leq
i<K\log N}\prod\limits_{\mu^{2}(d)=1,d|
ja^{i}+l, p|d\Rightarrow K<p}\frac{2^{\omega(d)}}{d})^{\frac{1}{2}}\\
&&\ll (K\log N\prod\limits_{p\in
P_{a,j,k,l}}(1-\frac{1}{p}))^{\frac{1}{2}}
(\sum\limits_{\mu^{2}(d)=1,p|d\Rightarrow K<p}\sum\limits_{1\leq i<K\log N, d|ja^{i}+l}\frac{2^{\omega(d)}}{d})^{\frac{1}{2}}\\
 &&\ll (K\log N\prod\limits_{p\in
P_{a,j,k,l}}(1-\frac{1}{p}))^{\frac{1}{2}}
(\sum\limits_{\mu^{2}(d)=1,p|d\Rightarrow K<p}\sum\limits_{1\leq
i<K\log N,
e_{a,j,k,l}(d)|i}\frac{2^{\omega(d)}}{d})^{\frac{1}{2}}\\
&&\ll (K\log N\prod\limits_{p\in
P_{a,j,k,l}}(1-\frac{1}{p}))^{\frac{1}{2}}
(\sum\limits_{\mu^{2}(d)=1,p|d\Rightarrow K<p}
  \frac{(K\log N)2^{\omega(d)}}{de_{a,j,k,l}(d)})^{\frac{1}{2}}\\
 &&\ll K\log N\prod\limits_{p\in
 P_{a,j,k,l}}(1-\frac{1}{p})^{\frac{1}{2}}
\end{eqnarray*}

So, we have \begin{eqnarray*}&&\sum\limits_{1\leq i\leq K\log
N,i\not\equiv I(a,j,k,l)\pmod {p}\,\,\mbox{for any}\,\,p\in
P_{a,j,k,l}}Q_{N,a,i,j,k,l}\\
&&\ll\frac{KN}{W\log N}\prod\limits_{3\leq p\leq
K}(1-\frac{2}{p})^{-1}\prod_{q|W}(1-\frac{2}{q})^{-1}
\prod\limits_{p\in P_{a,j,k,l}}(1-\frac{1}{p})^{\frac{1}{2}}.
\end{eqnarray*}

Thus, we get\begin{eqnarray*}&&\sum\limits_{1\leq i\leq K\log
N}Q_{N,a,i,j,k,l}\\
&&\ll\frac{KN}{W\log N}\prod\limits_{3\leq p\leq
K}(1-\frac{2}{p})^{-1}\prod_{q|W}(1-\frac{2}{q})^{-1}
\prod\limits_{p\in P_{a,j,k,l}}(1-\frac{1}{p})^{\frac{1}{2}}+\log N.
\end{eqnarray*}

 By Lemma 1, we get

\begin{eqnarray*}&&Q\geq c_{1}\frac{N}{W\log N}\prod_{q|W}(1-\frac{1}{q})^{-1}-c_{2}\frac{KN}{W\log
N}\prod\limits_{3\leq p\leq
K}(1-\frac{2}{p})^{-1}\prod_{q|W}(1-\frac{2}{q})^{-1}
\sum\limits_{a=2}^{K}\sum\limits_{(j,k,l)\in R}\prod\limits_{p\in P_{a,j,k,l}}(1-\frac{1}{p})^{\frac{1}{2}}\\
&&-\sum\limits_{(j,k,l)\in R}Q_{N,1,0,j,k,l}-O(\log N)\\
 &&\geq c_{1}\frac{N}{W\log N}\prod_{q|W}(1-\frac{1}{q})^{-1}-c_{3}\frac{KN}{W\log
N}\prod\limits_{3\leq p\leq
K}(1-\frac{2}{p})^{-1}\prod_{q|W}(1-\frac{1}{q})^{-2}
\sum\limits_{a=2}^{K}\sum\limits_{(j,k,l)\in R}\prod\limits_{p\in P_{a,j,k,l}}(1-\frac{1}{p})^{\frac{1}{2}}\\
&&\geq c_{1} \frac{N}{W\log
N}\prod_{q|W}(1-\frac{1}{q})^{-1}(1-c_{4}\prod\limits_{3\leq p\leq
K}(1-\frac{2}{p})^{-1}\prod_{q|W}(1-\frac{1}{q})^{-1}
\sum\limits_{a=2}^{K}\sum\limits_{(j,k,l)\in R}\prod\limits_{p\in P_{a,j,k,l}}(1-\frac{1}{p})^{\frac{1}{2}})\\
&&\geq c_{1}\frac{N}{W\log
N}\prod_{q|W}(1-\frac{1}{q})^{-1}(1-c_{5}(\log
K)^{2}\sum\limits_{a=2}^{K}\sum\limits_{(j,k,l)\in
R}\exp(\sum\limits_{q|W}\frac{1}{q}-\sum\limits_{p\in
P_{a,j,k,l}}\frac{1}{2p})\\
&&\geq c_{1}\frac{N}{W\log
N}\prod_{q|W}(1-\frac{1}{q})^{-1}(1-2c_{5}(\log
K)^{2}K^{4}\exp(K-\frac{M}{8K^{3}})).
\end{eqnarray*}

Taking $M>\max\{12K^{3}, 8K^{4}+8K^{3}\log (4c_{5}(\log
K)^{2}K^{4})\}$, we get $Q\geq C\frac{N}{W\log
N}\prod\limits_{q|W}(1-\frac{1}{q})^{-1}$, where the constant $C$ is
absolute.

This complets the proof of the Theorem 1.

\vskip 3mm{\bf Acknowledgement.}
   I am grateful to the referee for his/her useful suggestions on this paper.
   In October 2009, Professor Hao Pan in Nanjing gave a talk to
introduce his  work "On the number of distinct prime factors of
$nj+a^{h}k$". I would like to thank Professor Hao Pan for his
wonderful talk which attract my interest to this topic.


\begin {thebibliography}{30}



\bibitem{Chen6} Y. G. Chen and X. G. Sun,  \ { \it On Romanoff's
constant}, J. Number Theory 106(2004), 275-284.


\bibitem{Chenjnt5} Y. G. Chen, R. Feng and N. Templier,
 \ {\it Fermat numbers and integers of the form $a^{k}+a^{l}+p^{\alpha}$ },
Acta Arith. 135 (2008), 51-61.

\bibitem{CS}F. Cohen and J. L. Selfridge, \ {\it Note every number is the sum or difference of two
prime powers}, Math. Comput. 29(1975), 79-81.

\bibitem{Crocker} R. Crocker, \ {\it On the sum of a
prime and two powers of two}, Pacific J. Math. 36(1971), 103-107.

\bibitem{Erdos} P. Erd\H os,\ {\it On integers of the form $2^{r}+p$ and
some related problems}, Summa Brasil. Math. 2(1950), 113-123.

\bibitem{Halberstam} H. Halberstam and H. E. Richert,\ {\it
Sieve Methods}, Academic Press Inc. (London) Ltd., 1974.





\bibitem{MBN} M. B. Nathanson,  \ {\it Additive Number Theory, The Classical Bases},
Springer, New York, 1996.


\bibitem{pan} H. Pan, \ {\it On the integers not of the form $p+2^{\alpha}+2^{\beta}$},
Acta Arith. 148 (2011), no. 1, 55-61.

\bibitem{hpan} H. Pan, \ {\it On the number of distinct prime factors of  $nj+a^{h}k$},
Monatsh math. 175(2014), 293-305.

\bibitem{pintz} J. Pintz, \ {\it A note on Romanov's constant},
Acta Math. Hungar. 112(2006), 1-14.

\bibitem{Polignac} A. de Polignac, { \it Recherches nouvelles sur les nombres premiers},
 C. R. Acad. Sci. Paris Math.
29(1849), 397-401, 738-739.

\bibitem{Romanoff} N. P. Romanoff, {\it \"{U}ber einige S\"{a}tze der additiven
Zahlentheorie}, Math. Ann. 109(1934), 668-678.

\bibitem{RL} J. Barkley Rosser and Lowell Schoenfeld, {\it Approximate formulas for some
functions of prime numbers}, Illinois Journal of Mathematics,
6(1962), 64-94.


\bibitem{Sun}  Z. W. Sun and M. H. Le,
\ {\it Integers not of the form $c(2^{a}+ 2^{b}) + p^{\alpha}$},
Acta Arith. 99 (2001), 183-190.

\bibitem{Suntdai} X. G. Sun and L. X. Dai,
\ {\it Chen's conjecture and its generalization}, Chin. Ann. Math.
34B (2013), 957-962.

\bibitem{Tao}  T. Tao, \ {\it A remark on primality testing and decimal expansion}, J. Austr. Math.
Soc. \\3( 2011), 405-413.

\bibitem{Yuan} P. Z. Yuan, \ {\it Integers not of the form $c(2^{a}+2^{b})+p^{\alpha}$}, Acta Arith. 115
(2004), 23-28.


\end{thebibliography}
\end {document}